\documentclass[a4paper,12pt]{amsart}
\usepackage{amsmath}
\usepackage{amsmath}
\usepackage{amssymb}
\usepackage{amsthm}
\usepackage{graphicx,color,wrapfig}

\theoremstyle{plain}
\newtheorem{lemma}{Lemma}
\newtheorem{Definition}[lemma]{Definition}

\newtheorem{cor}[lemma]{Corollary}
\newtheorem{rem}[lemma]{Remark}

\newenvironment{Proof}[1][Proof]{
    \par
    \topsep 6pt plus 6pt
    \trivlist
    \item[\hskip\labelsep\bfseries #1.]\ignorespaces}%
    {\qed\endtrivlist
}

\newcommand{\R}{\mathbb{R}}

\newcommand{\Tr}{\text{Tr}}
\newcommand{\cov}{\text{cov}}

\begin{document}

\author{Hayk Mikayelyan, Henrik Shahgholian}
\title[Convexity of the free boundary]{Convexity of the free
boundary for an exterior free boundary problem involving the perimeter}



\thanks{The first author thanks G\"oran Gustafsson Foundation and ESF Programme on
``Global and Geometric Aspects of Non-Linear PDE'' for visiting
appointments to KTH, Stockholm}
\thanks{The second author is partially supported by Swedish Research Council}

\address{Hayk Mikayelyan \\Max-Planck-Institut f\"ur
Mathematik in den Naturwissenschaften\\
Inselstrasse 22\\
04103 Leipzig\\
Germany}
\email{hayk@mis.mpg.de}

\address{Henrik Shahgholian\\
Institutionen f\"or Matematik \\
Kungliga Tekniska H\"ogskolan\\
100 44 Stockholm \\
Sweden}
\email{henriksh@e.kth.se}

\begin{abstract}
We prove that if the given compact set $K$ is convex then
a minimizer of the functional
$$
I(v)=\int_{B_R} |\nabla v|^p dx+\text{Per}(\{v>0\}),\,1<p<\infty,
$$
over the set $\{v\in H^1_0(B_R)|\,\, v\equiv 1\,\,\text{on}\,\, K\subset B_R\}$
has a convex support, and as a result all its level sets are convex as well.
We derive the free boundary condition for the minimizers and prove that
the free boundary is analytic and the minimizer is unique.
\end{abstract}

\maketitle

\section{Introduction}

\subsection{The Problem}\label{firstsection}

The following problem has been considered in \cite{Maz}:
given a bounded domain $K\subset B_R\subset\R^n$ ($R$ large),
satisfying the interior ball condition,
find a (local) minimizer of the functional
\begin{equation}\label{funk}
I(v)=\int_{B_R} F(|\nabla v|)dx+\text{Per}(\{v>0\})
\end{equation}
over the set of functions
$\{v\in H^1_0(B_R)| v\equiv 1 \,\,\text{on}\,\,  K\}$,
where $F\in C^1([0,+\infty))$
is a positive convex function,
with $F(0)=0$ and for some $1<p<+\infty$ and $0<\lambda<\Lambda<+\infty$
$$
\lambda t^{p-1}\leq F'(t)\leq\Lambda t^{p-1}.
$$
Here we set
$\text{Per}(\{v>0\})=+\infty$ if $\chi_{\{v>0\}}\notin BV(\R^n)$.
This problem is the one-phase exterior
analogue of the problem introduced in \cite{ACKS}
for a functional with general convex function $F(t)$ in the first term
(in \cite{ACKS} they treat the case $F(t)=t^2$).

In \cite{ACKS} (two phase, $p=2$) and \cite{Maz} (one phase,
$1<p<\infty$) it is proved that the minimizers are Lipschitz
continuous. This gives that the free boundary
$\Gamma_u:=\partial\{x|u(x)>0\}$ is an almost minimal surface, thus
$C^{1,1/2}$-smooth. Let us recall some facts from the theory of
almost minimal surfaces following \cite{T}.

A set $\Omega$ has almost minimal boundary in $B_1$ if for every $A\Subset B_1$
there exist $R$, $0 < R < \text{dist}(A; \partial B_1)$ and
$$
\alpha:(0,R)\to[0,+\infty), \,\,\,\alpha(r)\downarrow_{r\to 0} 0,
$$
such that
$$
\text{Per}(\Omega;B_r(x))\leq\text{Per}(\Omega';B_r(x))+\alpha(r)r^{n-1}
$$
for every $x\in  A$, $r \in(0; R)$ and $\Omega'$ with $\Omega'\Delta\Omega\Subset B_r(x)$.

\begin{lemma} \cite{T}
Suppose $\Omega$ has almost minimal boundary in $B_1$ with $\alpha(r)=r^{2\lambda}$,
$\lambda\in(0,1/2]$. Then

(i) $\partial^*\Omega$ is a $C^{1,\lambda}$ hypersurface, and

(ii) $H^s(\partial\Omega\backslash\partial^*\Omega\cap B_1)=0$ for each $s>n-8$.
\end{lemma}Here $ \partial^* \Omega $ is the reduced boundary of $\Omega$
(see \cite{EG}).

As it is shown in \cite{ACKS} and \cite{Maz} (one phase,
$1<p<\infty$) the Lipschitz regularity of the minimizer gives that
the free boundary is an almost minimal surface with $\alpha(r) =r$
hence the reduced boundary $\Gamma_u^*:= \partial^* \{x|u(x)>0\}$
is $C^{1,\frac{1}{2}}$ regular and the singular set is of
Hausdorff dimension $n-8$ or less.

\begin{rem}\label{supplane}
Any blow-up at the almost minimal surface is a minimal cone
(see [T], p. 85).
Thats why all points of the free boundary at which we can
find a supporting smooth surface belong to the smooth part $\Gamma_u^*$.
\end{rem}

In this paper
we restrict ourselves to the case $F(t)=t^p$, $p>1$, i.e.,
the functional
\begin{equation}\label{funkp}
I(v)=\int_{B_R} |\nabla v|^p dx+\text{Per}(\{v>0\}),
\end{equation}
though we want to mention that the same ideas and methods
will work in the general case (\ref{funk})
if we put some additional (rather weak) conditions on the function $F$.

The main result of this paper is the following theorem.

\noindent \textbf{Theorem A}. \textsl{\ \ If $K\Subset B_R$ is a
convex set with non-empty interior and $u$ is a minimizer of
(\ref{funkp}) over the set $\{v\in H^1_0(B_R)| v\equiv 1
\,\,\text{on}\,\,  K\}$ then the minimizer is unique, the set
$\{u>0\}$ is convex and the free boundary $\partial\{u>0\}\cap B_R$
is an analytic surface.}

We also derive the free boundary condition in case of general
(non-convex) $K$ and prove that any minimizer $u_K$ satisfies the
following inclusion
$$\Omega_{u_K}\subset\Omega_{u_{\cov(K)}},$$ see the notations below.
This means for instance that for large $R$ we will have
$\Omega_{u_K}\Subset B_R$.

It is noteworthy that convexity results for the so-called
Bernoulli free boundary problem has been extensively studied, see
\cite{A}, \cite{HS1}, \cite{HS2} and the references therein.

\subsection{Notations}\label{Notat}
In the sequel we use the following notations:
\vspace{2mm}

\begin{tabbing}
$ \mathbf{R}_{+}^{n}$  \hspace{25mm}      \=  $\{x\in \mathbf{R}^{n}:x_{1}>0\}$\\
$  B(z,r) $  \>  $ \{x\in \mathbf{R}^{n}:|x-z|<r\} $, \\
$ B_{r} $  \>  $ B(0,r) $, \\
$ \chi _{D} $  \>  characteristic function of the set $ D $, \\
$ \partial D $  \> boundary of the set $ D $, \\
$ \Omega_u $  \>  $ \left\{ x\in \mathbf{R}^{n} :u\left( x\right)>0\right\} $, \\
$ \Gamma_u $  \>  $ \partial \Omega_u $ the free boundary,\\
$ \Gamma_u^* $  \>  $ \partial^* \Omega _u $ the reduced boundary of $\Omega_u$
(see \cite{EG}),\\
$ \cov(U) $  \> the convex hull of the set $ U$.
\end{tabbing}

\vspace{3mm}

\subsection{Organization of the paper}
In Sections \ref{energy} and \ref{sechopf} we develop some technical
tools which will be used in the proofs coming after. In Section
\ref{secfbc} we derive the free boundary condition and prove that
the reduced free boundary is analytic. An interesting geometric
result about the mean curvature of the boundary of the convex hull
of a non-convex domain is proved in Section \ref{secconcres}. The
main result of the paper, the convexity of the free boundary and the
uniqueness of the minimizer is proved in Section \ref{secmain}.

\section{An energy estimate for $p$-harmonic extensions}\label{energy}

Assume $K\Subset\Omega_1\subset\Omega_2$, where $K, \Omega_1,\Omega_2$
are open and bounded subsets of $\R^n$ with non-empty interior,
and that $u_j$ minimizes
the functional
\begin{equation}\label{Hfunk}
J(v)=\int |\nabla v|^p dx
\end{equation}
in the class of functions $\{v\in H^1_0(\Omega_j)|\,v\equiv
1\,\,\,\text{on}\,\,\, K\}$ ($j=1,2$). Then we say that $u_2$ is
the $p$-harmonic extension of $u_1$ from $\Omega_1$ to $\Omega_2$.
The use of word ``extension'' is a little bit misleading here,
since $u_1\not=u_2$ in $ \Omega_1$, but we keep it as an analogy
to extension by zero, which would be standard in this situation,
 because of zero boundary data on the boundary of $\Omega_1$.
We also extend functions in $H^1_0(\Omega)$ by zero and assume
that they are defined in all $\R^n$.

In this section we prove the following lemma.

\begin{lemma}\label{estimate}
If  $u_2$ is the $p$-harmonic extension of $u_1$ from $\Omega_1$ to
$\Omega_2$, where $\Omega_1$ and $\Omega_2$ have piecewise
$C^{1,\alpha}$ boundary. Then
\begin{multline}\label{norgnahat}
0\leq\int_{\Omega_2}|\nabla u_1|^p-|\nabla u_2|^pdx-
\int_{\Omega_2\backslash\Omega_1}(p-1)|\nabla u_2|^p dx\leq\\
-p\int_{\partial \Omega_1\backslash \partial \Omega_2}
u_2[|\nabla u_1|^{p-2}\partial_\nu u_1-|\nabla u_2|^{p-2}\partial_\nu u_2]dH^{n-1}.
\end{multline}
\end{lemma}
\begin{Proof}
We write
$$
v \Delta_p u =-p|\nabla u|^{p-2}\nabla u \nabla v +
p\text{div}(v |\nabla u|^{p-2}\nabla u),
$$
and using Gauss' theorem we obtain
\begin{multline}\label{Hharm1}
\int_{\Omega_2}p|\nabla u_2|^{p-2}\nabla u_2\nabla(u_1-u_2)dx=\\
-\int_{\Omega_2}(u_1-u_2)\Delta_p u_2dx+
\int_{\partial\Omega_2}p|\nabla u_2|^{p-2}(u_1-u_2)\partial_\nu u_2dH^{n-1}=0
\end{multline}
From here we have
\begin{multline*}
\int_{\Omega_2}|\nabla u_1|^p-|\nabla u_2|^p=\\
\int_{\Omega_2}(p-1)|\nabla u_2|^p +
[|\nabla u_1|^p-
p|\nabla u_2|^{p-2}\nabla u_1\nabla u_2  ]dx\\
=\int_{\Omega_2\backslash\Omega_1}(p-1)|\nabla u_2|^pdx+\\
\int_{\Omega_1}|\nabla u_1|^p-|\nabla u_2|^p-
p|\nabla u_2|^{p-2}\nabla u_2\nabla (u_1-u_2)dx.
\end{multline*}

Now we are going to estimate the last integral.
Let us consider the following function
$$
\Phi(t)=|\nabla(u_2+t(u_1-u_2))|^p.
$$
From the convexity and monotonicity of $t^p$ it follows that $\Phi$ is convex in $t$.
So we can write
$$
0\leq \Phi(1)-\Phi(0)-\Phi'(0)\leq \Phi'(1)-\Phi'(0).
$$
This gives us exactly the following
\begin{multline*}
0\leq |\nabla u_1|^p-|\nabla u_2|^p-
p|\nabla u_2|^{p-2}\nabla u_2\nabla (u_1-u_2)
\leq\\
 p|\nabla u_1|^{p-2}\nabla u_1\nabla (u_1-u_2)-
p|\nabla u_2|^{p-2}\nabla u_2\nabla (u_1-u_2)
\end{multline*}
in $\Omega_1$. Integrating partially in the domain $\Omega_1$ like
in (\ref{Hharm1}) we get that
\begin{multline*}
 \int_{\Omega_1}p|\nabla u_1|^{p-2}\nabla u_1\nabla (u_1-u_2)-
p|\nabla u_2|^{p-2}\nabla u_2\nabla (u_1-u_2)=\\
-p\int_{\partial \Omega_1\backslash \partial \Omega_2} u_2(|\nabla
u_1|^{p-2}\partial_\nu u_1-|\nabla u_2|^{p-2}\partial_\nu
u_2)dH^{n-1}
\end{multline*}
\end{Proof}

\begin{rem}\label{nonCalpha}
Note that the first inequality in (\ref{norgnahat}) does not require
smoothness assumptions on $\partial\Omega_{i}$, $i=1,2$.
\end{rem}

\section{The Hopf lemma for $p$-harmonic functions in domains with
Liapunov-Dini boundary}\label{sechopf}

Let us present the definition of Liapunov-Dini surface following
\cite{W}.

\begin{Definition}
A Liapunov-Dini surface $S$ is a closed, bounded
$(n-1)-$dimensional surface satisfying the following conditions:

(a) At every point of $S$ there is a uniquely defined tangent
(hyper-)plane, and thus also a normal.

(b) There exits a Dini modulus of continuity $\epsilon(t)$ such
that if $\beta$ is the angle between two normals, and $r$ is the
distance between their foot points, then the inequality
$\beta\leq\epsilon(r)$ holds.

(c) There is a constant $\rho>0$ such that for any point $x\in S$
any line parallel to the normal at $x$ meets $S\cap B_\rho(x)$ at
most once.
\end{Definition}
A modulus of continuity $\epsilon(r)\to 0$ as ${r\to 0}$ is called
Dini modulus of continuity if $\int t^{-1}\epsilon(t)dt<\infty$.

Note that a domain $E$ with Liapunov-Dini boundary satisfies a
kind of interior and exterior Dini condition in the following
sense: There exists a convex Liapunov-Dini domain $K$ such that
for any point $x_0\in \partial E$ there exists a translation and
rotation $K_{x_0}$ of the domain $K$ satisfying
$$
K_{x_0}\subset E, \,\,(K_{x_0}\subset\R^n\backslash
E)\,\,\,\text{and}\,\,\,\partial K_{x_0}\cap\partial E=\{x_0\}.
$$

The following lemma is proved in \cite{J}, see also \cite{MPS}.
\begin{lemma}\label{jmps}
Let $\Omega\backslash K$ be a convex ring and $u$
be its $p$-capacitary potential. Then
$$
\Delta_q u \le 0,\,\,\text{if}\,\,\, 1 < q \le p
$$
and
$$
\Delta_q u \ge 0, \,\,\text{if}\,\,\,  p \le q \le \infty.
$$
\end{lemma}
Now we formulate and prove the main result of this section, which
might be a known result but we could not find any reference.
\begin{lemma}\label{hopf-mopf}
Assume $u$ is a $p$-harmonic function in the domain $U$. Further
assume $y\in \partial U$, $ \partial U$ satisfies the interior and
exterior Dini conditions locally near $y$ and $u(x)\geq u(y)$ for
all $x\in U$. Then there exist positive constants $r_0,c,C$ such
that
$$
cr<\max_{\{x:|x-y|<r\}} u(x)-u(y)<Cr,
$$
for $0<r<r_0$.
\end{lemma}
\begin{Proof}
Let us take the function $w$ to be the minimizer of the Dirichlet
integral in $\{v\in H_0^1(K_2)|v\equiv 1
\,\,\,\text{on}\,\,K_1\}$, where $K_1$ and $K_2$ are convex
domains with Liapunov-Dini boundary and $K_1\Subset K_2$. Thus we
have $\Delta w=0$ on $K_2\backslash \overline{K_1}$. From the Hopf
lemma for harmonic functions (Thm. 2.5, \cite{W}) and the
convexity and regularity of the level sets of $w$ (see \cite{L})
we know that $\nabla w(x)\not= 0$, for any $x\in K_2\backslash
\overline{K_1}$. Now we will prove the existence of a smooth,
convex function $f: [0,1] \to [0,1]$, $f(0)=0$, $f(1)=1$ such that
$$
\Delta_p f(w)\geq 0
$$
in $K_2\backslash \overline{K_1}$ and $0<f'(t)<+\infty$ for all
$t\in[0,1]$. This will mean that the function $f(w)$ is a
sub-solution for $\Delta_p$ and has non-vanishing gradient, thus
it will work as a standard barrier function.

We have
\begin{equation*}
\Delta_p f(w)=p|\nabla f(w)|^{p-2}\Delta f(w) + p(p-2)|\nabla
f(w)|^{p-4}\Delta_\infty f(w),
\end{equation*}
where $$\Delta_\infty v=\sum_{i,j} v_{ij} v_i v_j$$ is the well
known infinity Laplace operator\footnote{In our definition by taking
$H(t)=t^p$ the operator $\Delta_p$ differs from the usual one by a
factor $p$.}. On the other hand
$$
\nabla f(w)=f'(w)\nabla w,
$$

$$
\Delta f(w)=f'(w)\Delta w+f''(w)|\nabla w|^2=f''(w)|\nabla w|^2,
$$

$$
\Delta_\infty f(w)=(f'(w))^3\Delta_\infty w+(f'(w))^2f''(w)|\nabla
w|^4.
$$

So we need to find a function $f $ such that
\begin{multline*}
\Delta_p f(w)=pf''(w)f'(w)^{p-2}|\nabla w|^p + \\
p(p-2)(f'(w)|\nabla w|)^{p-4}\left[ (f'(w))^3\Delta_\infty
w+(f'(w))^2f''(w)|\nabla w|^4\right]\geq 0,
\end{multline*}
or
\begin{equation}\label{bareranhav}
\frac{f''(w)}{f'(w)}\geq \frac{2-p}{p-1}|\nabla
w|^{-4}\Delta_\infty w.
\end{equation}
We see that for $p\geq 2$ we can take $f(t)\equiv t $. This
follows from the Lemma \ref{jmps}.

In case $1<p<2$ we continue as follows. We have from \cite{W} that
the derivatives of $w$ are continuous up to the boundary and do not
vanish. Moreover we have bounds for the second derivatives of $w$
near the boundary (formula (2.4.1) in \cite{W})
$$
|D^2w|\leq \zeta (d(x)),
$$
where $\zeta \in L^1(0,dist(K_1,\R^n\backslash K_2)/2)$, and
$d(x)$ is the distance function from the boundary of the domain
$K_2\backslash K_1$. Coming back to our case there exists a
function $\zeta_1(t)\in L^1((0,1))\cap C((0,1))$ such that
$|\nabla w|^{-4}|\Delta_\infty w|\leq\zeta_1(w)$ in $K_2\backslash
\overline{K_1}$. Let us now integrate (\ref{bareranhav}) in
$w\in[t,1]$,
$$
\int_t^1 \frac{f''(\tau)}{f'(\tau)}d\tau=
\int_{f'(t)}^{f'(1)}\frac{ds}{ s}\geq \frac{2-p}{p-1}\int_t^1
\zeta_1(\tau)d\tau.
$$
Thus we can take for instance
$$
f(t)=c \int_0^t \exp\left(-\frac{2-p}{p-1} \int_\tau^1
\zeta_1(s)ds\right)d\tau,
$$
where the constant $c>0$ is chosen to get $f(1)=1$.

Note that the function $1-f(w)$ is a super-solution for $\Delta_p$
in $K_2\backslash \overline{K_1}$ and will give us bounds from
above.
\end{Proof}

\begin{rem}\label{grad-mrad}
In case of the $C^{1,\alpha}$ boundary the existence of the
gradient of the function $u$ at the boundary is known (see
\cite{Li}) and we can write $$0<c<|\nabla u(y)|<C.$$
\end{rem}

\section{The free boundary condition}\label{secfbc}

First let us prove that for convex $K$ the free boundary stays away from the set
$K$.
\begin{lemma}\label{distfromK}Let the set $K$ be 
convex and $u$ be the minimizer of (\ref{funk}). Then
there exists a constant $\delta$ depending on $n$, $p$ and the set
$K$ such that $\text{dist}(x,K)\geq\delta$, for all $x\in\Gamma_u$.
\end{lemma}
\begin{Proof}
Let us take the points $y\in\Gamma_u$ and $x\in K$ such that
$\text{dist}(y,x)=r_0:=\text{dist}(\Gamma_u,K)$.

Let us denote by $V(r):=|B_r(x)\backslash \Omega_u|$ and
$A(r)=\mathcal{H}^{n-1}(\partial B_r(x)\backslash \Omega_u)$, so
that $V'(r)=A(r)$. Then we have from the isoperimetric inequality
and from the minimality condition that
$$
(V(r))^{\frac{n-1}{n}}\leq c \text{Per}(V(r))\leq 2c A(r).
$$
Now integrating
$$
2c\leq V'(r)(V(r))^{-\frac{n-1}{n}}
$$
in the interval $(r_0,r)$ we obtain
\begin{equation}\label{distanhav}
H^{n-1}(\partial B_r(x)\backslash \Omega_u)\geq c (r-r_0)^{n-1},
\end{equation}
where $c$ depends on the dimension.

From the convexity of $K$ and the fact that it has a non-empty
interior we know that there is a cone $\mathcal{C}$ and $r_K>0$ such
that for any point $y\in
\partial K$ there exists a rotation and translation of the set
$\mathcal{C}_{r_K}:= \mathcal{C}\cap B_{r_K}$ such that $0\mapsto y$
and $\mathcal{C}_{r_K}$ is mapped into $K$. In other words at any
point of $\partial K$ we can put a conical set of fixed opening and
length $r_K$ lying inside K. This gives that
$$
H^{n-1}(\partial B_r(y)\cap K )\geq c_K r^{n-1},
$$
for $r<r_K$ and all $y\in\partial K$, $c_K$ depends only on $K$. Let
us take the function $\zeta(x)$ to be the $p$-harmonic potential of
the convex ring $B_1(0)\backslash (B_{1/4}(0)\cap \mathcal{C} )$.

{\bf Step 1:} We first exclude the case $r_0=0$. Assume
$y\in\Gamma_u\cap K$. We take as a perturbation of $u$ the function
$v(x):=\max(u(x),\zeta_r(x))$, where $\zeta_r(x):=\zeta((x-y)/r)$
and we can without loss of generality assume that
$\{x|\zeta_r(x)=1\}\subset K$ for all $r<r_K$. Note that the
function $\zeta_r$ is the $p$-harmonic extension of the function $u$
from $\{u<\zeta_r\}\cap\Omega_u$ to $\{u<\zeta_r\}$, which together
with Lemma \ref{estimate} and Remark \ref{nonCalpha} gives the
following
\begin{multline}\label{estr_0}
(p-1)\int_{B_r(y)\backslash \Omega_u}|\nabla \zeta_r|^pdx\leq
\int_{\{u<\zeta_r\}}|\nabla u|^p-|\nabla \zeta_r|^p dx
\leq\\
H^{n-1}(\partial B_r(y)\backslash \Omega_u)-\text{Per}(\Omega_u;B_r(y)) \leq
C r^{n-1}.
\end{multline}
The second inequality uses the fact of $u$ being a minimizer. On the
other hand using (\ref{distanhav}) (remember that $r_0$ is assumed
to be 0) we obviously can find constants $c_1,c_2$ depending only on
$K$, $n$ and $p$ such that
\begin{multline}\label{xixi}
\int_{B_r(y)\backslash \Omega_u}|\nabla \zeta_r|^pdx\geq
\int_{B_r(y)\backslash (\Omega_u\cup B_{3r/4}(y))}|\nabla \zeta_r|^pdx\geq \\
c_1\int_{B_r(y)\backslash  B_{3r/4}(y)}|\nabla \zeta_r|^pdx\geq
c_2 r^{n-p},
\end{multline}
a contradiction. In the second and third inequalities of
(\ref{xixi}) we used (\ref{distanhav}) and the fact that
$cr^{-1}<|\nabla \xi_r|<Cr^{-1}$ in $B_r(y)\backslash  B_{3r/4}(y)$
for some positive constants depending on $r_K$.

{\bf Step 2:} Now we know that $r_0>0$ and we can use the estimates
used by Mazzone (Lemma 3.2, \cite{Maz}) for the terms in
(\ref{estr_0}). If we denote by $d'(x):=\text{dist}(x,\partial
B_r(y))$ we obtain that

\begin{multline}
H^{n-1}(\partial B_r(y)\backslash \Omega_u)-\text{Per}(\Omega_u;B_r(y)) \leq\\
-\int_{\partial^*(B_r(y)\backslash \Omega_u)}\langle\nabla d',\nu\rangle dH^{n-1}=\\
-\int_{B_r(y)\backslash \Omega_u}\Delta d'dx\leq \frac{c_n}{r_0}|B_r(y)\backslash \Omega_u|.
\end{multline}
On the other hand as in (\ref{xixi})
$$
\int_{B_r(y)\backslash \Omega_u}|\nabla \zeta_r|^pdx\geq
c |B_r(y)\backslash \Omega_u| r^{-p},
$$
where $c$ depends only on $n$ and $K$.
Summing up we obtain that
$$
c r_0^{-p}\leq \frac{c_n}{r_0},
$$
where all constants depend only on $n$, $p$ and $K$.
\end{Proof}

\begin{lemma}\label{lemfbc}
The reduced free boundary $\Gamma^*$ is analytic and $\Gamma^*\cap B_R$ 
satisfies the free boundary condition
\begin{equation}\label{fbhavas}
(p-1)|\nabla u|^p=\kappa(\Gamma^*_u),
\end{equation}
where $\kappa$ is the mean curvature.
Moreover on $\overline{\Omega}_u\cap\partial B_R$ we have pointwise
the inequality
\begin{equation}\label{fbanhavas}
(p-1)|\nabla u|^p\geq\kappa(\partial B_R).
\end{equation}
\end{lemma}
\begin{Proof}
{\bf Step 1:} We first derive the free boundary condition in the
weak sense using the domain variation method and show that
$\Gamma^*$ is $C^{2,\alpha}$ regular.

Assume the origin is a reduced free boundary point $0\in\Gamma^*$,
thus we can assume that for a small $\delta>0$ in the neighborhood
$\mathcal{N}_\delta:=\{x||x'|<\delta,\,|x_n|<\delta\}$ the free
boundary is a graph $\Gamma^*=\{x|x_n=\phi(x') \}$, with $\phi\in
C^{1,1/2}$ (see Section \ref{firstsection} and \cite{Maz}),
$\phi(0)=|\nabla\phi(0)|=0$.

For a vector field $\eta\in C^{1}_0(\mathcal{N}_\delta;\R^n)$, $\sup
\eta\le 1$ and small enough $\epsilon$ consider the bijective map
$\Phi_\epsilon(x)=x+\epsilon \eta(x)$ and the function
$u_\epsilon(y)=u(\Phi^{-1}_\epsilon(y))$.

From the minimality of
$u$ we have that
$$
\int_{B_{R}} |\nabla u_\epsilon(y)|^pdy-
\int_{B_{R}} |\nabla u(x)|^pdx + \text{Per}(\{u_\epsilon>0\})
-\text{Per}(\{u>0\})\geq 0.
$$
Let us now calculate the terms above. We are following the book of Ambrosio,
Fusco, Pallara (\cite{AFP}, page 360), where all these calculations are
carried out in a similar
situation.
Since
$$
\int_{B_{R}} |\nabla u_\epsilon(y)|^pdy= \int_{B_{R}} |\nabla
u(x)\cdot \nabla\Phi_\epsilon^{-1}(\Phi_\epsilon(x))|^p |\det
\nabla\Phi_\epsilon(x)|dx
$$
and
$$
\nabla\Phi_\epsilon^{-1}(\Phi_\epsilon(x))=I-\epsilon\nabla\eta(x)+o(\epsilon),
$$

$$
\det \nabla\Phi_\epsilon(x)=1+\epsilon\text{div}\eta(x)+o(\epsilon),
$$
we see that
\begin{multline*}
\int_{B_{R}} |\nabla u_\epsilon(y)|^pdy- \int_{B_{R}} |\nabla
u(x)|^pdx=\\ \epsilon \int_{B_R}\left( |\nabla u(x)|^p
\text{div}\eta(x)- p|\nabla u|^{p-2} \langle\nabla
u,\nabla\eta\cdot\nabla u\rangle \right)dx + o(\epsilon).
\end{multline*}
On the other hand
$$
\text{Per}(\{u_\epsilon>0\})
-\text{Per}(\{u>0\})=
\epsilon\int_{\Gamma^*}\text{div}^{\Gamma^*}\eta d\mathcal{H}^{n-1}+o(\epsilon),
$$
where $\text{div}^S F(x)=\sum_{k=1}^{n}\langle \nabla^S F_k(x), e_k \rangle$
is the tangential divergence of $F$ on surface $S$ and
$\nabla^S f$ is the projection of $\nabla f(x)$ on the tangent space $T_x S$
(see Definition 7.27 and Theorem 7.31 in \cite{AFP}).

Integrating by parts in $\mathcal{N}_\delta\cap\{u>0\}$ we obtain
\begin{multline*}
\int_{B_R} |\nabla u(x)|^p \text{div}\eta(x)dx=\\ \int_{\Gamma^*}
|\nabla u(x)|^p \langle \eta(x),\nu\rangle  d\mathcal{H}^{n-1} -
\int_{B_R} \langle \nabla|\nabla u(x)|^p ,\eta\rangle  dx,
\end{multline*}
where $\nu$ is the normal vector, and
\begin{multline*}
-\int_{B_{R}}
p|\nabla u|^{p-2} \langle\nabla u,\nabla\eta\cdot\nabla u\rangle dx=\\
\int_{B_{R}\cap\{u>0\}}\Delta_p u \langle\eta, \nabla u\rangle +
p|\nabla u|^{p-2}\langle\eta,\nabla^2 u \cdot  \nabla u\rangle dx- \\
\int_{\Gamma^*} p \langle\eta, \nabla u\rangle |\nabla u(x)|^{p-2}
\partial_\nu u d\mathcal{H}^{n-1}.
\end{multline*}
Noting that $\langle \nabla|\nabla u(x)|^p ,\eta\rangle = p|\nabla
u|^{p-2}\langle\eta,\nabla^2 \cdot u \nabla u\rangle $ and summing
up and letting $\epsilon$ go to 0 we obtain that
$$
\int_{\Gamma^*}\text{div}^{\Gamma^*}\eta d\mathcal{H}^{n-1}=
\int_{\Gamma^*} (p-1) |\nabla u(x)|^p \langle
\eta(x),\nu\rangle d\mathcal{H}^{n-1},
$$
for any $\eta\in C^1_0(\mathcal{N}_\delta;\R^n)$.

If we now rewrite the left hand side in terms of function $\phi$ and
use the Proposition 7.40 from \cite{AFP} we obtain that
\begin{equation}\label{curvelliptic}
-\text{div}\left(\frac{\nabla\phi(x')}{\sqrt{1+|\nabla\phi(x')|^2}}\right)=
(p-1)|\nabla u|^p (x',\phi(x'))
\end{equation}
weakly in $\{|x'|<\delta\}$. The $C^\alpha$ regularity of the right
hand side and the theory of quasilinear elliptic equations (see
\cite{GT}) give that $\phi$ is $C^{2,\alpha}$, so the reduced free
boundary is $C^{2,\alpha}$-smooth as well and the free boundary
condition (\ref{fbhavas}) is true pointwise on $\Gamma^*$.

{\bf Step 2:}
Now since higher regularity of the boundary implies the higher
regularity of the function $u$ up to the boundary, we can use the
bootstrapping argument and obtain arbitrary smoothness, so the
boundary is $C^\infty$.

{\bf Step 3:}
The analyticity follows from the theory of
elliptic coercive systems (see \cite{KNS}). Here we refer to
the paper of Argiolas (\cite{Ar}, p. 144), where a similar problem is treated
in all details.

{\bf Step 4:}
The inequality (\ref{fbanhavas}) is due to the fact that we can carry out the
domain variation only in one direction near $\partial B_R$.
\end{Proof}

\section{A concavity result}\label{secconcres}

From now on we denote by $\kappa(\partial U)$ the interior mean curvature
(in viscosity sense) of
the
$C^{1,1}$ part of the boundary of a domain $U$
as follows.
Assume $0\in \partial U$ and the interior normal $\nu_{\partial U}(0)$
shows in the
direction of the $e$-axis.
We take
$$
\kappa(\partial U)(0):=\inf_{\mathcal{A}\in\mathfrak{A}} \kappa
(S_\mathcal{A})(0),
$$
where $S_\mathcal{A}=\{(x,e)|e=\langle\mathcal{A} x,x \rangle\}$ and
$\mathfrak{A}$ is the set of all symmetric matrices $\mathcal{A}$
such that the set $S_\mathcal{A}$
(the graph of a quadratic polynomial) locally touches
$\partial U$ from inside.

Let us consider the convex hull $\cov(U)$ of a (non-convex) set $U$
with $C^{2}$ boundary. Note that then $\cov(U)$ has a
$C^{1,1}$ boundary (see \cite{KK}).
For notational reasons let us assume
$U\subset \R^{n+1}=\{(x,e)|x\in \R^n,e\in\R\}$.

The following lemma will be useful and is easy to
prove.

\begin{lemma}\label{lemkusc}
The function $\kappa(\partial \cov(U))(x)$ is
upper semi-continuous on $\partial \cov(U)$.
\end{lemma}

Assume we have a point $x_0\in\partial \cov(U)\backslash \partial U$, then
from the definition of the convex hull we know that $x_0$ is a convex combination of 
$n$ points from  $\partial \cov(U)\cap \partial U$, i.e., $x_0=\sum_{k=1}^n \alpha_k y_k$, 
$\sum_{k=1}^n\alpha_k=1$, $\alpha_k\geq 0$, $y_k\in\partial \cov(U)\cap \partial U$ 
for $k=1,\dots,n$.
Since $x_0\notin \partial U$ more than one of $\alpha_k$ will be different from zero,
thus there are points $y_0,z_0\in
\partial \cov(U)$ such that
$x_0$ lies in the interval $(y_0,z_0)\subset \partial \cov(U)$.

\begin{lemma}\label{lemmekbazhk}
The function
$$
\frac{1}{\kappa(\partial \cov(U))}(x)
$$ 
is concave on the interval $(y_0,z_0)\subset \partial \cov(U)$. 
Moreover if $\kappa(\partial \cov(U))(x)=0$
for some $x\in (y_0,z_0)$
then $\kappa(\partial \cov(U))(x)=0$ for all $x\in (y_0,z_0)$.
\end{lemma}

\begin{Proof}
We need to show that
$$
\frac{1}{\kappa(\partial \cov(U))}\Big(\frac{x^1+x^2}{2}\Big)\geq
\frac{1}{2}\Big( \frac{1}{\kappa(\partial \cov(U))}(x^1)
+\frac{1}{\kappa(\partial \cov(U))}(x^2) \Big)
$$
for all $x^1,x^2\in (y_0,z_0)$.
Without loss of generality we can assume
$x^1=(-1,0,\dots,0)$ and $x^2=(1,0,\dots,0)$.
Since the supporting planes of $\cov(U)$
at $x^1$ and $x^2$ coincide we can further assume that
the graphs of quadratic polynomials
$$
u=\langle\mathcal{A}_1 (x-x^1),(x-x^1) \rangle\,\,\, \text{and}\,\,\,
u=\langle\mathcal{A}_2 (x-x^2),(x-x^2) \rangle
$$
given by positive symmetric matrices
$\mathcal{A}_1$ and $\mathcal{A}_2$ locally
touch the boundary $\partial \cov(U)$ from inside and
$0<2\Tr\mathcal{A}_i-\kappa(x^i) <\epsilon$ for $i=1,2$.

\begin{figure}
\begin{center}
\input{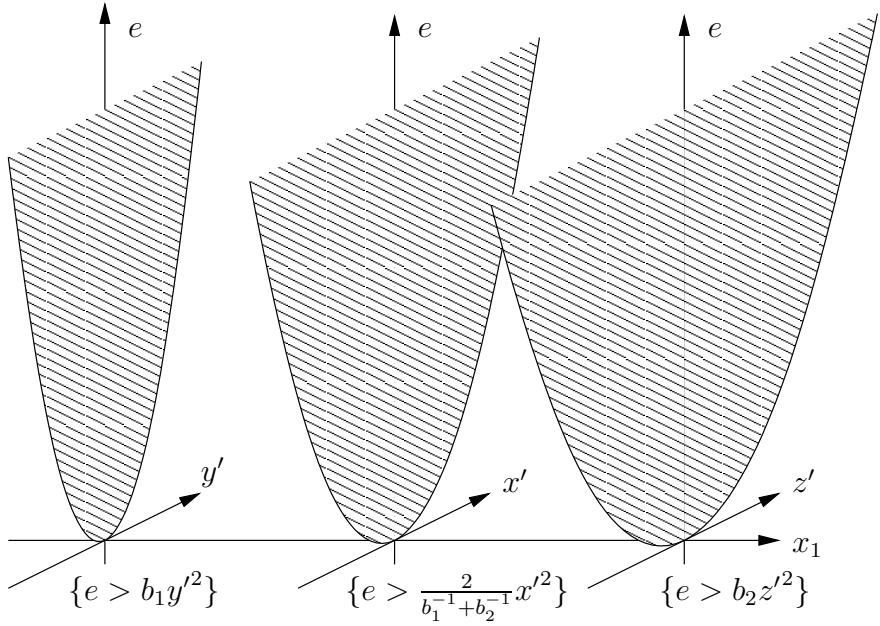}
\caption{Convex hull of two parabolas}
\label{parabolas}
\end{center}
\end{figure}

Since $x^1,x^2$ lie on the $x_1$ axis we can assume that for $i=1,2$
$$
\mathcal{A}_i=\left(\begin{array}{cccc}
a_i & 0 & \dots & 0 \\
0 & & & \\
\vdots & & \mathcal{B}_i & \\
0 & & &
\end{array}\right),
$$
where $\mathcal{B}_i $ are positive symmetric matrices and $0<a_i<\epsilon$.
The proof of this fact can be found in the Appendix.

Let us now consider the sets
$$
\{(- 1,x' , e)|
e>\langle\mathcal{B}_1 x',x' \rangle\}\,\,\,\text{and}\,\,\,
\{(1,x' , e)|
e>\langle\mathcal{B}_2 x',x' \rangle\},
$$
which touch the
boundary of $\cov(U)$ from inside locally at the points $x^1$
and $x^2$ respectively. Here $x'=(x_2,\dots,x_n)$.
We will now ``calculate'' the intersection
of the convex hull of this two sets with the plane $\{x|x_1=0\}$.
This will locally touch the boundary $\partial \cov(U)$ from inside
and give us the desired estimate on the mean curvature.
The intersection
of the convex hull of these two sets with the mentioned plane is
$\{(0,x',e)|e>u(x')\}$, where
\begin{equation}\label{infconv}
u(x')=\inf_{y'+z'=2x'}\frac{1}{2}\left(
\langle\mathcal{B}_1 y',y' \rangle
+\langle\mathcal{B}_2 z',z' \rangle
\right).
\end{equation}
We are going to calculate
explicitly the expression on the right hand side.
So for each $x'$ we are looking for the minimum of the following function
$$
w_{x'}(y')=\frac{1}{2}(
\langle\mathcal{B}_1 y',y' \rangle
+\langle\mathcal{B}_2 (2x'-y'),2x'-y' \rangle).
$$
After differentiation in $y'$ and some (simple) calculations we get
that the infimum in (\ref{infconv}) is attained at the values
$$
y'=2(\mathcal{B}_1+\mathcal{B}_2)^{-1}\mathcal{B}_2 x'
$$
and
$$
z'=2x'-y'=2(\mathcal{B}_1+\mathcal{B}_2)^{-1}\mathcal{B}_1 x'.
$$
Substituting now the values of $y'$ and $z'$ into (\ref{infconv})
and using the identity
$$
\mathcal{B}_1(\mathcal{B}_1+\mathcal{B}_2)^{-1}\mathcal{B}_2=
(\mathcal{B}_1^{-1}+\mathcal{B}_2^{-1})^{-1}
$$
we get
$$
u(x')=2
\langle(\mathcal{B}_1^{-1}+\mathcal{B}_2^{-1})^{-1} x',x' \rangle.
$$
Note that the invertibility of $\mathcal{B}_1+\mathcal{B}_2$
and $\mathcal{B}_1^{-1}+\mathcal{B}_2^{-1}$ follows from the
strict positivity of all eigenvalues of $\mathcal{B}_1,\mathcal{B}_2$.
In three dimensions, when matrices $\mathcal{B}_1,\mathcal{B}_2$
are given by positive numbers $b_1,b_2$, this interesting result is
illustrated in Figure \ref{parabolas}.

The proof now follows from the inequalities bellow:
\begin{multline}\label{3anhav}
\frac{2}{\kappa(\partial \cov(U))}\Big(\frac{x^1+x^2}{2}\Big)\geq
\frac{1}{2\Tr (\mathcal{B}_1^{-1}+\mathcal{B}_2^{-1})^{-1}}
\geq
\frac{1}{2\Tr \mathcal{B}_1} + \frac{1}{2\Tr \mathcal{B}_2}\\ \geq
 \frac{1}{\kappa(\partial \cov(U))(x^1)+\epsilon}
+\frac{1}{\kappa(\partial \cov(U))(x^2)+\epsilon}.
\end{multline}
 Note that
$\epsilon>0$ is arbitrary small and we have the first and the
third inequalities in
(\ref{3anhav}) by the construction of $\mathcal{B}_1,\mathcal{B}_2$
and from the properties of the convex hull.
The second inequality can be found in \cite{ALL}.

The case when  $\kappa(\partial \cov(U))(x)=0$
for some $x\in (y_0,z_0)$ follows from (\ref{3anhav}).
\end{Proof}

\section{Convexity of the free boundary}\label{secmain}

In the proof of the key Lemma \ref{lemkarevor} we will use the following lemma
(Lemma 4.1, \cite{LS}).
Let $K\subset U$
be a compact convex set, $U$ be open and non-convex
and $\cov(U)$ be the convex hull of $U$. Further assume
that the function $u$ minimizes the
functional (\ref{funkp}) over the set
$\{v\in H^1_0(\cov(U))| v\equiv 1\,\,\text{on}\,\, K\}$ and that
the segment $[y_0, z_0]\subset \partial\cov(U)$. Then the following
lemma is true.

\begin{lemma}\label{lemLS}
The function
$$
\frac{1}{|\nabla u|}(x)
$$ is convex on $(y_0, z_0)$.
\end{lemma}

This is due to the fact (see \cite{L}) that the level sets of a $p$-harmonic
potential in a convex ring are convex.

The following lemma is key to the proof of the main result.

\begin{lemma}\label{lemkarevor}
Let $u$ be a (local) minimizer of (\ref{funkp}) and denote by
$\cov(\Omega_u)$ the convex hull
of $\Omega_u$. Assume $u^c$ be the minimizer of
\begin{equation}\label{norfunk}
\int_{\cov(\Omega_u)\backslash K} |\nabla v(x)|^p dx
\end{equation}
over the set $\{v\in H^1_0(\cov(\Omega_u))|v\equiv 1\,\,\text{on}\,\,K\}$.

Then $\partial\cov(\Omega_u)$ is locally a $C^{1,1}$ surface and is a solution
of the (pointwise) free boundary inequality
\begin{equation}
(p-1)|\nabla u^c(x)|^p\geq\kappa(\partial\cov(\Omega_u)),
\end{equation}
where $\kappa$ is the interior mean curvature.
\end{lemma}
\begin{Proof}
The $C^{1,1}$ regularity of the $\partial\cov(\Omega_u)$ follows
from the fact that at all points of $\partial
\Omega_u\cap\partial\cov(\Omega_u)$ we have a supporting plane, thus
(Remark \ref{supplane}) $\partial \Omega_u$ is smooth in the
neighborhood of this points, i.e., all singular points of $\partial
\Omega_u$ have positive distance from $\partial\cov(\Omega_u)$. This
means that $\partial\cov(\Omega_u)$ is as regular as a convex hull
of a domain with smooth boundary, that is $C^{1,1}$ (see \cite{KK}).

We get the desired inequality on
$\partial\cov(\Omega_u)\cap\partial \Omega_u$
from the maximum principle and Lemma \ref{lemfbc}.

Assume now that $x_0\in\partial\cov(\Omega_u)\backslash\partial
\Omega_u$. From the definition of the convex hull it follows that we
can always write $x_0=\sum_{k=1}^m\alpha_k y_k$, $y_k\in
\partial\cov(\Omega_u)\cap\partial \Omega_u$, $\alpha_k>0$,
$\sum_{k=1}^m \alpha_k=1$, $2\leq m\leq n$.

We proceed by induction in $m$. Assume there exist
two points
$y_1,y_2\in \partial\cov(\Omega_u)\cap\partial \Omega_u$
such that
$y_1$, $x_0$ and $y_2$ lay on one line.

We need to show that
\begin{equation}\label{gcivra}
 \frac{1}{p-1}\left(\frac{1}{|\nabla u^c(x)|}\right)^p-
 \frac{1}{\kappa(\partial\cov(\Omega_u))(x_0)}\leq 0.
\end{equation}
We know that
$\frac{1}{|\nabla u^c(x)|}$ and thus
$\left(\frac{1}{|\nabla u^c(x)|}\right)^p$ is convex on $[y_1,y_2]$ (Lemma \ref{lemLS}).
Since  (\ref{gcivra}) is true at the
points $y_1$ and $y_2$
the proof follows from the concavity of
$\frac{1}{\kappa(\partial\cov(\Omega_u))(x)}$
on the line segment $(y_1,y_2)$ and its lower
semi-continuity (Lemmas \ref{lemkusc} and \ref{lemmekbazhk}).

The induction step $m\Rightarrow m+1$ finishes the proof.
\end{Proof}

\noindent \textbf{Theorem}. \textsl{\ \
If $K$ is convex and $u$ is a minimizer of (\ref{funkp}) then $\Omega_u$ is
also convex.}
\begin{Proof}
Suppose  $\Omega_u$ is not convex. Let us take $u^c$ and
$\cov(\Omega_u)$ as in Lemma \ref{lemkarevor} and assume $0\in
\text{int}K$. Further take $u_r^c(x):=u^c(rx) $,
$\cov(\Omega^r_u)=r^{-1}\cov(\Omega_u)$ and $r_0:=\inf
\{r>0|\cov(\Omega^r_u)\subset\Omega_u\}>1$. Assume
$\partial\cov(\Omega^{r_0}_u)$ touches $\partial\Omega_u$ at the
point $\tilde{x}$. First note that as in Remark \ref{supplane} we
have that $\tilde{x}$ is not on $\partial B_R$ and that
$\partial\Omega_u$ is analytic near $\tilde{x}$. We have now
\begin{multline}\label{hakas}
\kappa(\partial\cov(\Omega^{r_0}_u))(\tilde{x})\leq
r_0^{1-p}(p-1)|\nabla u_{r_0}^c|^p(\tilde{x})\leq \\
r_0^{1-p}(p-1)|\nabla u|^p(\tilde{x})=
r_0^{1-p}\kappa(\partial\Omega_u)(\tilde{x}),
\end{multline}
where the first inequality follows from Lemma \ref{lemkarevor}, the
second one from the comparison principle and the third equality is
the free boundary condition. On the other hand from the definition
of $r_0$ we get that
$\kappa(\partial\cov(\Omega^{r_0}_u))(\tilde{x})\geq
\kappa(\partial\Omega_u)(\tilde{x})$ and $r_0>1$, a contradiction.
\end{Proof}
\begin{cor}
The free boundary is an analytic
surface.
\end{cor}
\begin{cor}
Using the same method as in the proof of the theorem one can easily prove the
uniqueness of the minimizer by a contradiction argument. Note that in the
two-phase (interior) case (see \cite{ACKS}) the minimizer is not unique.
\end{cor}

\begin{cor}
Again by the same method one can prove that the domain $\Omega_u$
of a minimizer $u$ of the problem with
general compact set $K$ (even non-connected) is included in the
domain $\Omega_{\tilde u}$
of the minimizer ${\tilde u}$ of the problem with compact set $\cov(K)$.
\end{cor}

\begin{cor}
For large enough $R$
$$
\Omega_{u_K}\Subset B_R.
$$
\end{cor}
\begin{Proof}
Due to the previous corollary we need to prove this only for convex $K$.
If $y\in\partial B_R\cap\Gamma_{u_K}$ then by convexity the conical set
$C(y,K):=\{x|x\in[y,z], z\in K\}\subset\Omega_{u_K}$ and
$$
c R^2\leq\text{Per}C(y,K)\leq\text{Per}\Omega_{u_K},
$$
where the constant $c$ depends on the set $K$. This contradicts to the
fact that the total energy $I(u)$ should decrease with $R$.
\end{Proof}

\subsection*{Acknowledgement}
The first author is grateful to Prof. S. Luckhaus for valuable
discussions.

\section*{Appendix}

Assume $0\in L\subset \partial\Omega$, where $L$ is the segment connecting
$(-1,0,\dots,0)$ and $(1,0,\dots,0)$,
$\Omega$ is convex domain with $C^{1,1}$-boundary and
$\Omega\subset\{x\in\R^{n+1}|x_{n+1}>0\}$.
We will show here that
$$
\kappa(\partial U)(0):=\inf_{\mathcal{A}\in\mathfrak{A}} \kappa
(S_\mathcal{A})(0)=\inf_{\mathcal{A}\in\mathfrak{A}^*} \kappa
(S_\mathcal{A})(0)
$$
where $S_\mathcal{A}=\{(x,e)|e=\langle\mathcal{A} x,x \rangle\}$,
$\mathfrak{A}$ is the set of all symmetric matrices $\mathcal{A}$
such that the set $S_\mathcal{A}$
(the graph of a quadratic polynomial) locally touches
$\partial U$ from inside and $\mathfrak{A}^*\subset\mathfrak{A}$ is the subset
of the matrices of the form
$$
\mathcal{A}=\left(\begin{array}{cccc}
a & 0 & \dots & 0 \\
0 & & & \\
\vdots & & \mathcal{B} & \\
0 & & &
\end{array}\right).
$$
Note that the boundary of the set $\Omega$ around $0$ can be locally
given as graph of the following function
$$
u=f_{x_1}(x')+o(|x'|^2), \,\,\text{as}\,\,x'\to 0
$$
where $x'=(x_2,\dots,x_n)$ and
$f_{x_1}$ are homogeneous functions of order two, i.e.,
$f_{x_1}(x')=f_{x_1}(\frac{x'}{|x'|})|x'|^2$.
It is enough to show that the interior mean curvature of $\partial \Omega$
and of $\{x\in \R^{n+1}|x_{n+1}=f_{0}(x')\}$
at point $0$ is the same.
To see this let as fix any unit vector
$e=(\alpha_{1},\dots,\alpha_{n})$ in $\R^n$ and note that
\begin{multline*}
f_{\alpha_{1}t}(\alpha_{2}t,\dots,\alpha_{n}t)-
f_{0}(\alpha_{2}t,\dots,\alpha_{n}t)=\\
(f_{\alpha_{1}t}(\alpha_{2},\dots,\alpha_{n})-
f_{0}(\alpha_{2},\dots,\alpha_{n}))t^2=o(t^2),
\end{multline*}
as $t\to 0+$.



\begin{thebibliography}{A}


\bibitem[A]{A} A. Acker \textit{On the existence of convex
classical solutions for multilayer free boundary problems with general
nonlinear joining conditions} Trans. Amer. Math. Soc.  350  (1998),  no.
8, 2981-3020

\bibitem[Al]{Al} F. J., Jr. Almgren \textit{Existence and regularity
almost everywhere of
solutions to elliptic variational
problems with constraints} Mem. Amer. Math. Soc. 4 (1976), no. 165


\bibitem[ALL]{ALL} O. Alvarez, J.-M. Lasry, P.-L. Lions
\textit{Convex viscosity solutions and state constraints}
J. Math. Pures Appl., 76, 1997, 265-288


\bibitem[AFP]{AFP} L. Ambrosio, N. Fusco, D. Pallara, Diego
\textit{Functions of bounded variation and free discontinuity problems}
Oxford Mathematical Monographs. Oxford University Press, New York, 2000


\bibitem[Ar]{Ar} R. Argiolas
\textit{A two-phase variational problem with curvature}
Matematiche (Catania) 58 (2003), no. 1, 131--148


\bibitem[ACKS]{ACKS} I. Athanasopoulos, L. A. Caffarelli,
C. Kenig, S. Salsa \textit{An area-Dirichlet integral
minimization problem} Comm. Pure Appl. Math. 54 (2001), no. 4, 479-499

\bibitem[EG]{EG}  L. C. Evans, R. F. Gariepy \textit{Measure
theory and fine properties
of functions} Studies in Advanced Mathematics. CRC Press, Boca Raton, FL, 1992

\bibitem[GT]{GT} Gilbarg, David and Trudinger, Neil S.\textit{
Elliptic partial differential equations of second order.}
Springer Verlag, Berlin, 2001

\bibitem[HS1]{HS1} A. Henrot, H. Shahgholian \textit{Existence of
classical solutions to a free boundary problem for the
$p$-Laplace operator. I. The exterior convex case} J. Reine Angew. Math.
 521  (2000), 85--97

\bibitem[HS2]{HS2} A. Henrot, H. Shahgholian \textit{The one
phase free boundary problem for the $p$-Laplacian with non-constant
Bernoulli boundary condition} Trans. Amer. Math. Soc.  354  (2002),  no.
6, 2399--2416

\bibitem[J]{J} U. Janfalk \textit{Behaviour in the limit,
as $p\to\infty$, of minimizers of functionals involving
$p$-Dirichlet integrals} SIAM J. Math. Anal. 27(2) (1996), 341-360

\bibitem[KNS]{KNS} D. Kinderlehrer, L. Nirenberg, J. Spruck
\textit{Regularity in elliptic free boundary problems I}
J. Analyse Math. 34 (1978), 86--119

\bibitem[KK]{KK} B. Kirchheim, J. Kristensen \textit{Differentiability of convex
envelopes} C. R. Acad. Sci. Paris Sér. I Math.  333  (2001),  no. 8, 725--728

\bibitem[LS]{LS} P. Laurence, E. Stredulinsky \textit{Existence of regular
solutions with convex levels for semilinear elliptic equations with nonmonotone
$L^1$ nonlinearities. Part I} Indiana Univ. Math. J. vol. 39 (4) (1990),
1081-1114

\bibitem[L]{L} J. L. Lewis
\textit{Capacitary functions in convex rings}
Arch. Rational Mech. Anal.  66  (1977), no. 3, 201--224

\bibitem[Li]{Li} G. Lieberman
\textit{Boundary regularity for solutions of degenerate elliptic equations}
Nonlinear Anal.  12  (1988),  no. 11, 1203--1219

\bibitem[MPS]{MPS} J. Manfredi, A. Petrosyan, H. Shahgholian
\textit{A free boundary problem for $\infty$-Laplace equation}
Calc. Var. Partial Differential Equations 14 (2002), no. 3, 359--384

\bibitem[Maz]{Maz} F. Mazzone \textit{A single phase variational problem
involving the area of level surfaces} Comm. Part. Diff. Eq.
Vol. 28 (2003), no. 5\&6, 991-1004

\bibitem[T]{T} I. Tamanini \textit{Regularity results for almost minimal oriented
hypersurface in $R^n$} Quaderni del
Dipartimento di Matematica, Università di Lecce 1 (1994).

\bibitem[W]{W} K.-O. Widman \textit{Inequalities for the Green function
and boundary continuity of the gradient of solutions of elliptic differential
equations} Math. Scand. 21, 1967, 17-37 (1968)

\end{thebibliography}
\end{document}